\documentclass{article}
\usepackage{amsthm}
\usepackage{amsfonts,amssymb,amsmath}
\usepackage{enumitem}
\usepackage{titlesec}
\usepackage{fouriernc}
\usepackage[T1]{fontenc}
\usepackage{chngcntr}
\counterwithin{figure}{section}

\newlist{gcases}{enumerate}{1}
\setlist[gcases,1]{
  label={{\it Case}~{\it \Alph*}.},
  topsep=0ex,
  leftmargin=0in,
  labelsep=.1in,
  itemindent=.7in,
  itemsep=0ex 
}
\newlist{tenumerate}{enumerate}{1}
\setlist[tenumerate,1]{
  label={(\arabic*)},
  topsep=0ex,
  leftmargin=.3in,
  labelsep=.1in,
  itemindent=0in,
  itemsep=0ex
}

\newlist{titemize}{enumerate}{1}
\setlist[titemize,1]{
  label={$\bullet$},
  topsep=0ex,
  leftmargin=.3in,
  labelsep=.1in,
  itemindent=0in,
  itemsep=0ex
}

\titleformat{\section}
  {\normalfont\bfseries}   % The style of the section title
  {}                             % a prefix
  {0pt}                          % How much space exists between the prefix and the title
  {Section \thesection\quad}    % How the section is represented

\titleformat{name=\section,numberless}
  {\normalfont\bfseries}
  {}
  {0pt}
  {}

\newlength{\tabwidth}
\newlength{\tabheight}
\newlength{\tabrule}
\newlength{\tabwidthx}
\newlength{\tabheightx}

\def\gentabbox#1#2#3#4{\vbox to \tabheight{\setlength{\tabrule}{#3}%
  \setlength{\tabwidthx}{#1\tabwidth}\addtolength{\tabwidthx}{\tabrule}%

\setlength{\tabheightx}{#2\tabheight}\addtolength{\tabheightx}{-\tabheight}%
  \hbox to #1\tabwidth{%
 \hspace{-0.5\tabrule}\rule{\tabrule}{#2\tabheight}\hspace{-\tabrule}%
    \vbox to #2\tabheight{\hsize=\tabwidthx%
      \vspace{-0.5\tabrule}\hrule width\tabwidthx height\tabrule%
      \vspace{-0.5\tabrule}\vfil%
      \hbox to \tabwidthx{\hss#4\hss}%
        \vfil\vspace{-0.5\tabrule}%
      \hrule width\tabwidthx height\tabrule\vspace{-0.5\tabrule}}%
 \hspace{-\tabrule}\rule{\tabrule}{#2\tabheight}\hspace{-0.5\tabrule}}%
  \vspace{-\tabheightx}}}
\def\genblankbox#1#2{\vbox to \tabheight{\vfil\hbox to
#1\tabwidth{\hfil}}}
\def\tabbox#1#2#3{\gentabbox{#1}{#2}{0.4pt}{\strut #3}}

\catcode`\:=13 \catcode`\.=13 \catcode`\;=13
\catcode`\>=13 \catcode`\^=13
\def:#1\\{\hbox{$#1$}}
\def.#1{\tabbox{1}{1}{$#1$}}
\def>#1{\tabbox{2}{1}{$#1$}}
\def^#1{\tabbox{1}{2}{$#1$}}
\def;{\genblankbox{1}{1}\relax}
\catcode`\:=12 \catcode`\.=12 \catcode`\;=12
\catcode`\>=12 \catcode`\^=7

\newcommand\T{\mathbf{T}}
\newcommand\U{\mathbf{U}}

\begin{document}

\noindent
{\bf \large Isotypic components of left cells in type $D$}

\vspace{.3in}
\noindent
WILLIAM~M.~MCGOVERN \\
\vspace{.1in}
{\it \small Department of Mathematics, University of Washington, Seattle, WA, 98195}

\section*{Introduction}
This paper extends the results of \cite{M20} to type $D$, showing that the rule given there for computing basis vectors of isotypic components of a Kazhdan-Lusztig left (or right) cell extends to that type.

\section{Operators on tableau pairs}
We will continue to use the notation of the three parts of Garfinkle's series of papers on the classification of primitive ideals in types $B$ and $C$, as maintained and supplemented in \cite{MP16}.  Let $W$ be a Weyl group of type $D_n$.  If $\alpha,\beta$ are adjacent simple roots, now necessarily of the same length, then we use the wall-crossing operator $T_{\alpha\beta}$ as defined in \cite[2.1.10]{G92}, unless $\{\alpha,\beta\} = \{\alpha_1',\alpha_3\}= \{e_2+e_1,e_3-e_2\}$, in which case we use the operators $T_{\alpha_1'\alpha_3}^L,T_{\alpha_3,\alpha_1'}^L$ defined in \cite[4.3.5]{MP16}, alternatively denoting either of these by $T_{\alpha\beta}^L$.  These last operators do not preserve right tableaux in general; we will use them only in the cases where they do preserve this tableau, that is, either case (1) of \cite[4.3.5]{MP16}, or in case (2) of this definition or its analogue in case (4) if the cycle of the 3-domino is closed in  $\tilde\T_1$, or in case (3) or its analogue in case (4) if the cycle of the 3-domino is closed in $\T_1$.   As in \cite{M20}, we extend the operators $T_{\alpha\beta}$ to tableau pairs $(\T_1,\T_2)$ by decreeing that they act trivially on $\T_2$; the operators $T_{\alpha_1'\alpha_3}^L,T_{\alpha_3\alpha_1'}^L$ are already defined on tableau pairs.  We also use the operators $T_{\alpha X}^L,T_{X\alpha}^L$ defined in \cite[4.4.10]{MP16} and their truncations $U_{\alpha X}^L,U_{X\alpha}^L$, where the value or values of the latter operators on a tableau pair $(\T_1,\T_2)$ consist of all pairs in $T_{\alpha X}^L(\T_1,\T_2)$ or $T_{X\alpha}^L(\T_1,\T_2)$ whose right tableau is $\T_2$.  As in \cite{M20}, we need a family of additional operators to generate the equivalence relation among tableau pairs of having the same right tableau.  By analogy with the corresponding definition in \cite{M20}, we define a tableau shape to be a {\sl quasi staircase} if it takes the form $\lambda_n=(2n+1,\ldots,n+3,n+2,n+2,n,n,n-1,\ldots,1)$ or $\mu_n=(2n+2,2n+1,\ldots,n+3,n+2,n+2,n,n,n-1,\ldots,1)$; note that the quasi staircase shapes in type $D$ are just the transposes of the quasi staircase shapes in type $C$.   We then define operators $S_n,S_n',{}^tS_n$, and ${}^tS_n'$ much as we did in \cite{M20}, for each $n$ fixing two domino tableaux $\tilde\T_1,\tilde\U_1$ of shape $\lambda_n$ and two others $\tilde\T_2,\tilde\U_2$ of shape $\mu_n$, such that the two largest dominos of all of them lie vertically at the end of the two rows of length $n$ and horizontally at the end of the row just above these rows and each $\tilde\U_i$ is obtained from $\tilde\T_i$ by interchanging these dominos.  Then the operator $T_n$ is defined on a tableau pair $(\T_1,\T_2)$ such that the first $(n+1)^2$ dominos form the subtableau $\tilde\T_1$ or $\tilde\U_1$ by interchanging the two largest dominos in the subtableau while leaving all other dominos in $\T_1$ and $\T_2$ unchanged; we define $S_n'$ similarly, working with $\tilde\T_2$ and $\tilde\U_2$ rather than $\tilde\T_1$ and $\tilde\U_1$, and we define ${}^tS_n,{}^tS_n'$ to be the transposes of $S_n,S_n'$, respectively, taking the transpose of a pair on which $\S_n$ or $\S_n'$ is defined to the transpose of the image of that pair under $S_n$ or $S_n'$.  The composition $T_\Sigma$ of a sequence $\Sigma$ of operators $T_{\alpha\beta},T_{\alpha\beta}^L,U_{\alpha X}^L,U_{X\alpha}^L,S_n,S_n',{}^tS_n$, and ${}^tS_n'$, is defined as in \cite{M20}.

\section{Transitivity of the action on tableau pairs}

We now extend Theorem 1 of \cite{M20} to type $D$.

\newtheorem*{thm1}{Theorem 1}
\begin{thm1}
Given two pairs $(\T_1,\T_2),(\T_1',\T_2')$ of domino tableaux of the same shape such that $\T_2=\T_2'$, there is a sequence $\Sigma$ of operators $T_{\alpha\beta},T_{\alpha\beta}^L,U_{X\alpha}^L,U_{\alpha X}^L,S_n,\break S_n',{}^tS_n$, and ${}^tS_n'$, such that $(\T_1',\T_2')$ is one of the pairs in $T_\Sigma(\T_1,\T_2)$.
\end{thm1}

\begin{proof}
This is proved in the same way as in \cite{H97} and \cite[Theorem 1]{M20}, using Lemma 4.6.8 and Theorem 4.6.2 of \cite{MP16} in place of Theorem 3.2.2 of \cite{G93}.  Once again the most difficult case in proving the analogue of Lemma 4.6.8 is case I and it is this case that gives rise to the operators $S_n,S_n'$ and their transposes.
\end{proof}

We enlarge the operators $S_n,S_n',{}^tS_n$, and ${}^tS_n'$ to operators $T_n,T_n',{}^tT_n$, and $ {}^tT_n'$ as in \cite{M20}, so that these last operators are defined on all tableau pairs $(\T_1,\T_2)$ such that $\T_1$ can be moved through open cycles to produce a tableau with quasi staircase or transposed quasi staircase shape and the image of a tableau pair under an operator is either one or two tableau pairs of the same shape.  These operators again extend to $W$-equivariant linear maps from left cells on which they are defined to other left cells; the maps also preserve right cells.  Similarly the operators $U_{X\alpha}^L,U_{\alpha X}^L$ extend to $W$-equivariant maps $T_{X\alpha}^L,T_{\alpha X}^L$ from left cells of $W$ to other left cells that also preserve right cells.

\section{Decomposition of left cells}
The rule stated before Lemma 1 of \cite{M99} for constructing bases of isotypic components of Kazhdan-Lusztig left cells in types $B$ or $C$ continues to hold for a Weyl group $W$ of type $D$.  Fix left cells $\mathcal C,\mathcal R$ of $W$ lying in the same double cell $\mathcal D$.  Let $x$ be the unique element of $\mathcal C\cap\mathcal R$ whose left tableau $T_L(x)$ has special shape.  Let $\sigma$ be either of the two partitions of $2n$ corresponding to a representation $\pi$ of $W$ occurring in both $\mathcal C$ and $\mathcal R$, if $\sigma$ is not very even; if it is very even, then $\mathcal C$ and $\mathcal R$ are both irreducible as $W$-modules and there is nothing to do.  Let $e_1,\ldots,e_r$ be the extended open cycles of $T_L(x)$ relative to $T_R(x)$ such that moving $T_L(x)$ through these open cycles produces a tableau of shape $\sigma$.  Given any $w\in\mathcal C\cap\mathcal R$, let $T_L(w)$ be obtained from $T_L(x)$ by moving through the extended open cycles  $f_1,\ldots,f_s$ (relative to $T_R(x)$).  Put $\sigma_w=\pm1$ according as an even or odd number of $f_i$ appear among the $e_j$.  Set $R_\sigma = \sum_{w\in\mathcal C\cap\mathcal R} \sigma_wC_w$.

\newtheorem*{thm2}{Theorem 2}
\begin{thm2}
The right or left $W$-submodule generated by $R_\sigma$ is irreducible and $W$ acts on it by $\pi$.
\end{thm2}

\begin{proof}
This is proved in the same way as \cite[Theorem 2]{M20}, using \cite[4.4.10]{MP16} in place of \cite[2.3.4]{G92}, to show that the formula for $R_\sigma$ is compatible with the operators $T_{\alpha\beta},T_{\alpha\beta}^L,T_{\alpha X}^L,T_{X\alpha}^L$, in the sense that if a particular $R_\sigma$ coincides with $R_\sigma'$ and so transforms by $\pi$, then the same will be true of the image of $R_\sigma$ under any composition $\Sigma$ of maps $T_{\alpha\beta},T_{\alpha\beta}^L,T_{\alpha X}^L,T_{X\alpha}^L$ that is nonzero on $R_\sigma$.  We must show that this continues to hold for compositions including the operators $T_n,T_n',{}^tT_n$, and ${}^tT_n'$, and for this it is enough to check that the formula for $R_\sigma$ holds for particular cell intersections arising in the definition of the $S_n$ and $S_n'$, as in \cite{M20}, assuming inductively that it holds in smaller rank. 

Assume now that $W$ is of type $D_9$ and consider a cell intersection $I=\mathcal C\cap\mathcal C^{-1}$, where $\mathcal C$ is represented by an element with left tableau $\T_1$ chosen as above for the shape $(5,4,4,2,2,1)$; suppose in addition that the $7$-domino in $\T_1$ is horizontal and lies directly above the $8$-domino, while the $6$-domino is horizontal and lies at the end of the first row, so that the irreducible constituents of $\mathcal C$ as a $W$-module are indexed by the partitions $(5,5,3,3,1,1),(5,4,4,2,1,1),(5,5,3,2,2,1)$, and $(5,4,4,3,1,1)$ of 18.  Denote by $x_p$ the unique element in the intersection $I$ whose left and right tableaux have shape $p$.  Then compositions of operators $T_{\alpha\beta},T_{\alpha\beta}^L,U_{X\alpha}^L$, and $U_{\alpha X}^L$ act transitively on tableau pairs with a fixed right tableau not of shape $(5,4,4,2,2,1)$ or its transpose $(6,5,3,3,1)$, so the argument in the proof of \cite[Theorem 1]{M99} applies to show that

\begin{align*}
R_{5,5,3,3,1,1)}' = x_{(5,5,3,3,1,1)}+x_{(5,4,4,2,2,1)}+x_{(5,5,3,2,2,1)}+x_{(5,4,4,3,1,1)}\\ R_{(5,5,3,2,2,1)}' = x_{(5,5,3,3,1,1)}-x_{(5,4,4,2,1,1)}-x_{(5,5,3,2,2,1)} + x_{(5,4,4,3,1, 1)}
\end{align*}
\vskip .2in
\noindent while either

\begin{align*}
R_{(5,4,4,2,2,1)}' = x_{(5,5,3,3,1,1)}+x_{(5,4,4,2,2,1)}-x_{(5,5,3,2,2,1)}-x_{(5,4,4,3,1,1)}\\
 R_{(5,4,4,3,1,1)}' = x_{(5,5,3,3,1,1)}-x_{(5,4,4,2,2,1)}+x_{(5,5,3,2,2,1)}-x_{(5,4,4,3,1,1)}
 \end{align*}
\vskip .2in
\noindent or else 

\begin{align*}
R_{(5,4,4,2,2,1)}' = x_{(5,5,3,3,1,1)}-x_{(5,4,4,2,2,1)}+x_{(5,5,3,2,2,1)}-x_{(5,4,4,3,1,1)}\\
R_{(5,4,4,3,1,1)}' = x_{(5,5,3,3,1,1)}+x_{(5,4,4,2,2,1)}-x_{(5,5,3,2,2,1)}-x_{(5,4,4,3,1,1)}
\end{align*}
\vskip .2in
\noindent But now if we take the basis elements for the Weyl group $W'$ of type $D_5$ corresponding to the tableau pairs consisting of the first five dominos of every tableau in all the pairs corresponding to elements of $I$ and label the resulting elements $y_q$ in type $D_5$ by partitions $q$ of 10 as we did the elements of $I$ by partitions of 18, we find that $y_{(3,3,1,1,1,1)} - y_{(3,2,2,1,1,1)}$ transforms by the representation corresponding to $(3,2,2,1,1,1)$ of $W'$, whose truncated induction to $W$ is the direct sum of the representations corresponding to $(5,4,4,3,1,1)$ and $(5,4,4,2,2,1)$.  In order to make $x_{(5,4,4,2,2,1)} - x_{(5,5,3,3,1,1)}$ transform by representations lying in this last truncated induced representation, we find that the first pair of equations for $R_{(5,4,4,2,2,1)}' $ and $R_(5,4,4,3,1,1)'$ must hold, as desired..  It follows that the operators $T_2,T_2',{}^tT_2$, and ${}^tT_2'$, extended to linear maps between left cells regarded as $W$-modules, are indeed equivariant for the left $W$-action.  Like the maps $T_{\alpha\beta}$ and $T_{\alpha\beta}^L$, they are compatible with the formula for $R_\sigma$.  Similar arguments show that the linear extensions of the other maps $T_n,T_n',{}^tT_n$, and ${}^tT_n'$ are also $W$-equivariant and compatible with the formula for $R_\sigma$.  As in \cite{M20} we now have enough $W$-equivariant maps between left cells to validate the proof of \cite[Theorem 1]{M99}.
\end{proof} 

We then correct the statements of Theorem 4.2 in \cite{M96} and Theorems 2.1 and 2.2 in \cite{M'96} as in \cite{M20}, all three of these results being corrected and superseded by the following one:  {\sl given any two left cells $\mathcal C_1,\mathcal C_2$ in a Weyl group $W$ of type $D$ that have a representation $\pi$ of $W$ in common, there is a composition $\Sigma$ of maps $T_{\alpha\beta},T_{\alpha\beta}^L,T_{X\alpha}^L,T_{\alpha X}^L,T_n,T_n',{}^tT_n$, and ${}^tT_n'$ from $\mathcal C_1$ to $\mathcal C_2$ whose restriction to the copy of $\pi$ in $\mathcal C_1$ maps it isomorphically onto the corresponding copy of $\pi$ in $\mathcal C_2$}.


\begin{thebibliography}{10}


\bibitem{G90}
D.~Garfinkle.
\newblock On the classification of primitive ideals for complex classical {L}ie
  algebras (I),
\newblock {\em Compositio Math.}, 75(2):135--169, 1990.

\bibitem{G92}
D.~Garfinkle.
\newblock On the classification of primitive ideals for complex classical {L}ie
  algebras (II),
\newblock {\em Compositio Math.}, 81(3):307--336, 1992.

\bibitem{G93}
D.~Garfinkle.
\newblock On the classification of primitive ideals for complex classical Lie
  algebras (III),
\newblock {\em Compositio Math.}, 88:187--234, 1993.

\bibitem{H97}
B.~Hopkins.
\newblock Domino Tableaux and Single-Valued Wall-Crossing Operators,
\newblock Ph.D. dissertation, University of Washington, 1997.

\bibitem{L87}
G.~Lusztig.
\newblock Leading coefficients of character values of Hecke algebras,
\newblock {\em Proc.\ Symp.\ Pure Math.}, 47 (2): 235--262, 1987.

\bibitem{M96} 
W.~M.~McGovern.
\newblock Left cells and domino tableaux in classical Weyl groups,
\newblock {\em Compositio Math.}, 101:77--98, 1996.

\bibitem{M'96}
W.~M.~McGovern,
\newblock Standard domino tableaux and asymptotic Hecke algebras,
\newblock {\em Compositio Math.}, 101:99--108, 1996.

\bibitem{M99}
W.~M.~McGovern,
\newblock Errata and a new result on signs,
\newblock {\em Compositio Math.}, 117:117--121, 1999.

\bibitem{M20}
W.~M.~McGovern,
\newblock A family of operators generating domino tableaux of a fixed shape and a decomposition of left cells into isotypic components,
\newblock preprint, 2020.


\bibitem{MP16}
W.~M.~McGovern and T.~Pietraho,
\newblock On the classification of primitive ideals for complex classical Lie algebras (IV),
\newblock preprint, 2016.




\end{thebibliography}
\end{document}